\documentstyle[amssymb]{amsart}

\title{Bqo is \boldmath\P12-complete}
\author{Alberto Marcone}
\address{Dipartimento di Matematica, Universit\`a di Torino,
via Carlo Alberto 10, 10123 Torino, Italy}
\email{marcone@@dm.unito.it}

\newtheorem{theorem}{Theorem}[section]
\newtheorem{lemma}[theorem]{Lemma}
\newtheorem{corollary}[theorem]{Corollary}
\newtheorem{sublemma}{Sublemma}[theorem]
\theoremstyle{definition}
\newtheorem{definition}[theorem]{Definition}

\renewcommand{\P}[2]{$\Pi_{#2}^{#1}$}
\renewcommand{\S}[2]{$\Sigma_{#2}^{#1}$}

\newcommand{\Sq}[1]{\left[{#1}\right]^{\omega}}
\newcommand{\sq}[1]{\left[{#1}\right]^{<\omega}}
\newcommand{\Bai}{\cal N}
\newcommand{\Can}{\cal C}
\newcommand{\set}[2]{\left\{\,{#1}\mid{#2}\,\right\}}
\newcommand{\imply}{\;\longrightarrow\;}
\newcommand{\sse}{\longleftrightarrow}
\newcommand{\lh}{\operatorname{lh}}
\newcommand{\dom}{\operatorname{dom}}
\newcommand{\base}{\operatorname{base}}
\newcommand{\lt}{\left\langle}
\newcommand{\rt}{\right\rangle}
\newcommand{\conc}{{}^\smallfrown}
\newcommand{\sconc}{\conc\!\!}
\newcommand{\N}{\Bbb N}
\newcommand{\init}{\sqsubset}
\newcommand{\initeq}{\sqsubseteq}
\newcommand{\m}[1]{#1^-}
\newcommand{\R}{\,R\,}
\newcommand{\T}{T^*(C)}
\newcommand{\db}{\,\underline\ll\,}
\newcommand{\tri}{\vartriangleleft}
\newcommand{\ntri}{\ntriangleleft}
\newcommand{\f}{\tilde f}
\newcommand{\s}{\sigma}
\newcommand{\C}{D_x}
\newcommand{\CF}{C_{x,y}}
\newcommand{\Cf}{C^*_{x,y}}

\begin{document}
\pagestyle{plain}
\maketitle

\begin{abstract}
In this paper we give a proof of the \P12-completeness of the set of
countable bqos (viewed as a subset of the Cantor space). This result was
conjectured by Clote in \cite{Clote} and proved by the author in his Ph.d.\
thesis (\cite{thesis}, see also \cite{fbqo}): here we prove it using
Simpson's definition of bqo (\cite{SGS}) and as little bqo theory as
possible.
\end{abstract}

\section{Introduction}

Let $\Can$ be the Cantor space $2^\omega$ of the infinite sequences of $0$'s
and $1$'s with the product topology (where $2 = \{0,1\}$ is given the
discrete topology). Similarly $\Bai$ is the Baire space $\N^\omega$ with the
product topology (again $\N$ is endowed with the discrete topology). $\Can$
and $\Bai$ are Polish spaces (complete separable metric spaces). The basic
reference for this spaces from the viewpoint of descriptive set theory is
{\sc Moschovakis}' monograph (\cite{Mosch}).

\begin{definition}
A subset of a Polish space $\cal X$ is \S11 if it is the projection of a
Borel subset of $\cal X \times \cal Y$ where $\cal Y$ is a Polish space, it
is \P1n if it is the complement of a \S1n set and it is \S1{n+1} if it is the
projection of a \P1n subset of $\cal X \times \cal Y$.
\end{definition}

\begin{definition}
For $n \geq 1$ we say that a set $C \subseteq \cal X$ is {\em
\P1n-complete\/} if it is \P1n and for every \P1n set $P \subseteq \Bai$
there exists a continuous function $f: \Bai \to \cal X$ such that for every
$x \in \Bai$ we have $x \in P \sse f(x) \in C$. $f$ is called a {\em
reduction\/} of $P$ to $C$.
\end{definition}

A \P1n-complete set is true \P1n, i.e.\ it is \P1n but not \S1n. For an
excellent survey on true \P1n and \S1n sets see \cite{Becker}.

In this paper we will give an example of a subset of the Cantor space $\Can$
which is \P12-complete. This set arises from combinatorics and is
particularly interesting because it is a generalization of the canonical
example of a \P11-complete set, the set of all countable well-orderings. For
other examples of \P12-complete sets occurring in different branches of
mathematics see \cite{Becker} and the references quoted there.

A quasi-ordering (hereafter qo) consists of a set together with a binary
relation on its elements which is reflexive and transitive: this definition
is slightly more general than that of partial ordering (which requires also
anti-symmetry) and is useful whenever there is no canonical way of choosing
a representative among the elements that are equivalent under the given
relation. The most natural generalization of the concept of well-ordering to
qos is the notion of well quasi-ordering (hereafter wqo): a qo is a wqo if
it is well founded and contains no infinite set of mutually incomparable
elements. It is immediate that the notion of wqo is \P11; moreover the linear
orders that are wqos are exactly the well-orderings and hence the notion of
wqo is \P11-complete.

The concept of wqo is very natural but it does not enjoy nice closure
properties: this was discovered in the 1950's and in the 1960's {\sc
Nash-Williams} (\cite{NW1}, \cite{NW2}) proposed the stronger, but less
natural, notion of better quasi-ordering (hereafter bqo). Since then bqos
have become an interesting topic of research and very often a qo is proved
to be a wqo by showing that it is actually a bqo (one of the most famous
results of this kind is {\sc Laver}'s proof of Fra\"\i ss\'e's conjecture,
see \cite{Laver} and \cite{SGS}). For surveys of wqo and bqo theory see
\cite{Milner} and \cite{Pou-sur}. An alternative definition of bqo,
equivalent to {\sc Nash-Williams}' original one but without its combinatorial
flavor, has been given by {\sc Simpson} (\cite{SGS}) and has proved to be
very useful. Here, in contrast with \cite{thesis} and \cite{fbqo}, we will
use {\sc Simpson}'s definition because its more descriptive set theoretic
flavor may be more appealing to the intended readers of this paper. We
postpone the rather technical definition of bqo to section~\ref{sec:bqo}.

An immediate Tarski-Kuratowski computation shows that the set of all
countable bqos (viewed as a subset of $\Can$) is \P12 and a natural
conjecture made by {\sc Clote} (\cite{Clote}) stated that it is indeed
\P12-complete. We proved this conjecture in \cite{thesis} and \cite{fbqo}
using techniques originally devised to answer some questions dealing with the
fine analysis of the notion of bqo. Feeling that the extraction from
\cite{thesis} or \cite{fbqo} of all details of the proof of the
\P12-completeness of bqo might be arduous for the reader interested more in
descriptive set theory than in bqo theory, in this paper we give an
exposition of the proof which uses the least possible amount of bqo theory.

In section~\ref{sec:bqo} we give {\sc Simpson}'s definition of bqo and prove
some basic results of bqo theory which will be needed in the proof of the
main theorem. In particular, using {\sc Simpson}'s definition, we give a
proof of {\sc Pouzet}'s theorem (theorem~\ref{th:pouzet}) which allows us,
as far as bqo theory is concerned, to substitute arbitrary binary relations
in place of qos. In section~\ref{sec:smooth} we introduce the notion of
smooth subset of $\sq\N$ and prove some basic facts about it.
Section~\ref{sec:main} contains the proof of the main result of the paper,
i.e.\ theorem~\ref{main theorem}.

I am very much indebted to Stephen G. Simpson, who introduced me to bqo
theory, supervised my Ph.d.\ thesis where these results were originally
proved and suggested several improvements to a previous version of this
paper. I wish also to thank Maurice Pouzet, who explained to me his theorem
which is an important tool in the proof of the main result of this paper, and
Alessandro Andretta, whose interest in this result stimulated the writing of
this paper.

\section{Better quasi-orderings}\label{sec:bqo}

If $s$ is a finite sequence we denote by $\lh(s)$ its length and, for every
$i< \lh(s)$, by $s(i)$ its ($i+1$)-th element. We also write $s = \lt s(0),
\ldots, s(\lh(s) -1)\rt$, so that $\lt \rt$ denotes the empty sequence. If
$s$ and $t$ are finite sequences we write $s \initeq t$ if $s$ is an initial
segment of $t$, i.e.\ if $\lh(s) \leq \lh(t)$ and $\forall i < \lh(s) \;
s(i)=t(i)$. $s \init t$ has the obvious meaning and we extend this notation
also to the case where $t$ is an infinite sequence.

We write $s \conc t$ for the concatenation of $s$ and $t$, i.e.\ the sequence
$u$ such that $\lh(u) = \lh(s) + \lh(t)$, for every $i < \lh(s)$ $u(i) =
s(i)$ and for every $i < \lh(t)$ $u(\lh(s) + i) = t(i)$. If $s$ is a finite
sequence and $i \leq \lh(s)$ we denote by $s[i]$ the initial segment of $s$
of length $i$, i.e.\ the unique sequence $t$ such that $t \initeq s$ and
$\lh(t) = i$.

If $A$ is an infinite set we denote by $\Sq A$ the set of all countable
infinite subsets of $A$, by $\sq A$ the set of all finite subsets of $A$ and
by $\left[A\right]^n$ the set of elements of $\sq A$ with $n$ elements. We
identify members of $\sq \N$ and $\Sq \N$ with the sequences  (finite or
infinite) which enumerate them in increasing order. With this identification
$\Sq \N$ can be viewed as a closed subspace of $\Bai$ and is actually
homeomorphic to $\Bai$ via the map which sends $\alpha \in \Bai$ to $\set{k+
\sum_{i=0}^k \alpha(i)}{k \in \N} \in \Sq \N$. For any $A \in \Sq \N$ we will
always consider $\Sq A$ endowed with the topology arising from this
identification. A basis for this topology is given by the collection of all
sets of the form $N_s = \set{X \in \Sq A}{s \init X}$ where $s \in \sq A$.

\begin{definition}
If $X \subseteq \N$ is nonempty we denote by $\m X$ the set obtained from $X$
by removing its least element.
\end{definition}

The map $X \mapsto \m X$ is continuous from $\Sq \N$ in itself.

\begin{definition}
Let $(Q, \preceq)$ be a quasi-ordering and equip $Q$ with the discrete
topology. A {\em $Q$-array\/} is a Borel measurable function $f: \Sq A \to
Q$, where $A \in \Sq{\N}$. $f$ is {\em good (with respect to $\preceq$)\/}
if there exists $X \in \Sq A$ such that $f(X) \preceq f(\m X)$. If $f$ is not
good we say that it is {\em bad\/}. $f$ is {\em perfect\/} if for every $X
\in \Sq A$ we have $f(X) \preceq f(\m X)$.
\end{definition}

\begin{definition}
$(Q, \preceq)$ is {\em bqo\/} if every $Q$-array is good.

In this definition no role is played by the fact that $\preceq$ is a qo and
we can replace it by any binary relation on $Q$, which we usually denote by
$R$; in this case we say that $R$ is a {\em better binary relation\/} or {\em
bbr\/}.
\end{definition}

In the sequel if $f: \Sq A \to Q$ is a $Q$-array and $B$ is an infinite
subset of $A$ we will call the restriction of $f$ to $B$ the $Q$-array that
should be more precisely called the restriction of $f$ to $\Sq B$.

We will need the following classical result of descriptive set theory, known
as the Galvin-Prikry theorem (\cite{GP}): for its proof see e.g.\ \cite{SGS}.

\begin{theorem}\label{th:GP}
Let $A \in \Sq\N$ and suppose $\cal B$ is a Borel subset of $\Sq A$. Then
there exists $B \in \Sq A$ such that either $\Sq B \subseteq \cal B$ or $\Sq
B \cap \cal B = \emptyset$.
\end{theorem}

\begin{corollary}\label{bad/perfect}
Let $f: \Sq A \to Q$ be a $Q$-array and $R$ be a binary relation on $Q$. Then
there exists $B \in \Sq A$ such that $f$ restricted to $B$ is either bad or
perfect with respect to $R$.
\end{corollary}
\begin{pf}
Let $\cal B = \set{X \in \Sq A}{f(X) \R f(\m X)}$. Since $f$ is Borel
measurable and $X \mapsto \m X$ is continuous $\cal B$ is a Borel set. By
theorem~\ref{th:GP} there exists $B \in \Sq A$ such that either $\Sq B
\subseteq \cal B$ or $\Sq B \cap \cal B = \emptyset$: in the first case $f$
restricted to $B$ is perfect with respect to $R$, in the second case it is
bad with respect to $R$.
\end{pf}

\begin{lemma}\label{intersection}
If $R$ and $R'$ are two bbrs on the same set $Q$ then the relation $S= R \cap
R'$ is also bbr.
\end{lemma}
\begin{pf}
Let $f$ be a $Q$-array with $\dom(f) = \Sq A$: since $R$ is bbr by
corollary~\ref{bad/perfect} $f$ restricted to some $B \in \Sq A$ is perfect
with respect to $R$. Since $R'$ is bbr there exists $X \in \Sq B$ such that
$f(X) \,R'\, f(\m X)$ and hence $f(X) \,S\, f(\m X)$.
\end{pf}

A consequence of the Galvin-Prikry theorem is the following result, which
follows from the results of section~6 of \cite{Mathias} and is proved in
\cite{SGS}.

\begin{theorem}\label{th:Mathias}
Let $A \in \Sq\N$ and suppose $Y$ is a metric space and $f: \Sq A \to Y$ is
a Borel measurable function. Then there exists $B \in \Sq A$ such that the
restriction of $f$ to $B$ is continuous.
\end{theorem}

\begin{corollary}\label{cont=borel}
A binary relation $R$ on a set $Q$ is bbr if and only if every continuous
$Q$-array is good with respect to $R$.
\end{corollary}
\begin{pf}
One direction of the equivalence is trivial. For the other, if $f: \Sq A \to
Q$ is a $Q$-array by theorem~\ref{th:Mathias} (since $Q$ with the discrete
topology is metrizable) there exists $B \in \Sq A$ such that the restriction
of $f$ to $B$ is continuous. By hypothesis $f$ restricted to $B$ is good and
hence $f$ is good.
\end{pf}

One of the basic tools for showing that a qo is bqo is a theorem known as the
minimal bad array lemma. It is implicit in {\sc Nash-Williams}' work and its
present formulation (in terms of Borel measurable functions) was given by
{\sc Simpson} (\cite{SGS}). The proof we give here is due to {\sc van
Engelen}, {\sc Miller} and {\sc Steel} (\cite{vEMS}). Before proving the
minimal bad array lemma we prove a combinatorial lemma that is needed in its
proof and make two preliminary definitions which are needed for its
statement.

\begin{lemma}\label{combinatorial}
Let $\omega_1$ denote the first uncountable ordinal and
$\set{A_\alpha}{\alpha < \omega_1}$ be a sequence of elements of $\Sq \N$
such that $\alpha < \beta < \omega_1$ implies that $A_\beta \setminus
A_\alpha$ is finite. Then there exist $B \in \Sq \N$ and $I \in
\Sq{\omega_1}$ such that $B \subseteq \bigcap_{\alpha \in I} A_\alpha$.
\end{lemma}
\begin{pf}
Notice that the hypothesis of the lemma implies that if $s \in \sq{\omega_1}$
and $\beta = \max s$ then
$$\bigcap_{\alpha \in s} A_\alpha = A_\beta \setminus \bigcup_{\alpha \in s
\setminus \{\beta\}} \left(A_\beta \setminus A_\alpha\right)$$
is infinite.

We define by induction $I_n \in \left[ \omega_1 \right]^n$, $B_n \in \left[
\N \right]^n$ and $V_n \subseteq \omega_1$ uncountable such that
$$I_n \subset I_{n+1} \text{, } B_n \subset B_{n+1} \text{, } B_n \subseteq
\bigcap_{\alpha \in I_n} A_\alpha \text{, } B_n \subseteq \bigcap_{\alpha \in
V_n} A_\alpha \text{ and } I_n \cap V_n = \emptyset \text{.}$$
To complete the proof it will then suffice to let $I = \bigcup_n I_n$ and $B
= \bigcup_n B_n$.

We start by setting $I_0 = B_0 = \emptyset$ and $V_0 = \omega_1$. Supposing
$I_n$, $B_n$ and $V_n$ have already been defined let $\alpha_n$ be the least
element of $V_n$ and set $I_{n+1} = I_n \cup \{\alpha_n\}$. For every $\beta
\in V_n$ let $C_\beta = A_\beta \cap \bigcap_{\alpha \in I_{n+1}} A_\alpha$.
By the remark at the beginning of the proof $C_\beta$ is infinite and hence
there exists $m_\beta \in C_\beta \setminus B_n$. Since $V_n$ is uncountable
there exist $m \in \N$ and an uncountable $V_{n+1} \subseteq V_n \setminus
\{\alpha_n\}$ such that for every $\beta \in V_{n+1}$ we have $m_\beta = m$.
Let $B_{n+1} = B_n \cup \{m\}$. It is easy to check that $I_{n+1}$, $B_{n+1}$
and $V_{n+1}$ satisfy all the conditions.
\end{pf}

\begin{definition}
Let $(Q,\preceq')$ be a qo: if $f$ and $g$ are $Q$-arrays with domains $\Sq
A$ and $\Sq B$ respectively, we write $f \preceq' g$ if $B \subseteq A$ and
for all $X \in \Sq B$ we have that $g(X) \preceq' f(X)$. We write $f \prec'
g$ if $B \subseteq A$ and for all $X \in \Sq B$ we have that $g(X) \prec'
f(X)$ (i.e.\ $g(X) \preceq' f(X)$ and not $f(X) \preceq' g(X)$).
\end{definition}

\begin{definition}
Let $(Q,\preceq)$ be a qo: another qo $\preceq'$ on $Q$ is {\em compatible
with $\preceq$\/} if it is well founded and $q_0 \preceq' q_1$ implies $q_0
\preceq q_1$. In this setting we say that the $Q$-array $f$ is {\em minimal
bad\/} if it is bad (with respect to $\preceq$) and every $Q$-array $g$
satisfying $g \prec' f$ is good.
\end{definition}

Now we can state and prove the minimal bad array lemma.

\begin{theorem}\label{mba}
Let $\preceq$ and $\preceq'$ be qos on $Q$ such that $\preceq'$ is compatible
with $\preceq$. If $f$ is a $Q$-array which is bad (with respect to
$\preceq$) then there exists a minimal bad $Q$-array $g$ such that $g
\preceq' f$.
\end{theorem}
\begin{pf}
Suppose the theorem fails: we define by induction a sequence of bad
$Q$-arrays $\set{f_\alpha: \Sq{A_\alpha} \to Q}{\alpha < \omega_1}$ such that
if $\alpha < \beta < \omega_1$ then $f_\alpha \prec' f$, $A_\beta \setminus
A_\alpha$ is finite and for every $X \in \Sq{A_\alpha \cap A_\beta}$ we have
$f_\beta(X) \prec' f_\alpha(X)$. An application of lemma~\ref{combinatorial}
then easily shows that $\preceq'$ is not well founded, a contradiction.

Start by letting $f_0 = f$. The definition of $f_\alpha$ is immediate when
$\alpha$ is a successor ordinal: if $\alpha = \beta +1$ it suffices to take
$f_\alpha \prec' f_\beta$ bad.

If $\alpha < \omega_1$ is a limit ordinal and $f_\beta$ has been defined for
all $\beta < \alpha$ we begin by constructing a bad $Q$-array $g$ as follows.
Since $\alpha$ is countable we can reorder $\set{A_\beta}{\beta < \alpha}$
as $\set{B_n}{n \in \N}$ with $B_0 = A_0$. As in the proof of
lemma~\ref{combinatorial} we have that for every $k \in \N$ the set
$\bigcap_{n \leq k} B_n$ is infinite. Let $m_0$ be the least element of $B_0$
and $m_{k+1}$ be the least element of $\bigcap_{n \leq k+1} B_n$ which is
larger than $m_k$. Then $A = \set{m_k}{k \in \N}$ satisfies $A \subseteq A_0$
and $A \setminus A_\beta$ finite for every $\beta < \alpha$. We define $g:
\Sq A \to Q$ by $g(X) = f_\beta (X)$ where $\beta$ is maximal such that $X
\in \Sq{A_\beta}$ (notice that for every $X \in \Sq \N$ the set $\set{\beta
< \alpha}{X \in \Sq{A_\beta}}$ is finite because $\preceq'$ is well founded).

\begin{sublemma}
$g$ is a bad $Q$-array.
\end{sublemma}
\begin{pf}
To see that $g$ is Borel measurable it suffices to recall that $Q$ is given
the discrete topology and notice that for every $q \in Q$ and $X \in \Sq A$
we have $g(X) = q$ if and only if
$$\exists \beta < \alpha (X \in \Sq{A_\beta} \land X \in f_\beta^{-1}(\{q\})
\land \forall \beta' (\beta < \beta' < \alpha \imply X \notin
\Sq{A_{\beta'}}))$$
and hence $g^{-1}(\{q\})$ is a Borel set in $\Sq A$.

To see that $g$ is bad suppose that for some $X \in \Sq A$ we have $g(X)
\preceq g(\m X)$: if $g(X) = f_\beta(X)$ and $g(\m X) = f_{\beta'}(\m X)$,
from $\m X \subset X$ follows $\beta \leq \beta'$ and therefore
$f_{\beta'}(\m X) \preceq' f_\beta(\m X)$. Since $\preceq'$ is compatible
with $\preceq$ we have $f_{\beta'}(\m X) \preceq f_\beta(\m X)$ and hence
$f_\beta(X) \preceq f_\beta(\m X)$, contradicting the badness of $f_\beta$.
\end{pf}

Applying the successor step to $g$ we obtain $f_\alpha: \Sq{A_\alpha} \to Q$
which is bad and such that $f_\alpha \prec' g$. Since for every $\beta <
\alpha$ and $X \in \Sq{A \cap A_\beta}$ we have either $g(X) = f_\beta(X)$
or $g(X) \prec' f_\beta(X)$ it follows that $f_\alpha(X) \prec' f_\beta(X)$
for every $X \in \Sq{A_\alpha \cap A_\beta}$ and in particular $f_\alpha
\prec' f$. This completes our construction and the proof of the theorem.
\end{pf}

The following result allows us, whenever we are looking for a bqo, to look
instead for a bbr knowing that inside it we will find a bqo. It is a
consequence of a sharper theorem due to {\sc Pouzet} (\cite{Pou-pap}, see
also \cite{thesis} or \cite{fbqo}): here we translate {\sc Pouzet}'s proof
in the terminology of {\sc Simpson}'s definition of bqo.

\begin{theorem}\label{th:pouzet}
Let $R$ be a binary relation on a countable set $Q$. Then there exists a
partial ordering $\preceq$ of $Q$ such that $\preceq \; \subseteq R$ (as
subsets of $Q \times Q$) and such that $\preceq$ is bqo if and only if $R$
is bbr.

If we fix an enumeration of $Q$ we can view $R$ and $\preceq$ as elements of
$\Can \times \Can \cong \Can$: then the function $R \mapsto \preceq$ is
continuous.
\end{theorem}
\begin{pf}
Let $\{q_n\}$ be an enumeration of $Q$ and let $\#:Q \to \N$ be the function
such that $\#(q_n) = n$. We define $\preceq$ by primitive recursion as
follows: if $m \leq n$ suppose we have already established whether $q_i
\preceq q_j$ for all $i<m$ and $j \leq n$ and set $q_m \preceq q_n$ if and
only if $q_m \R q_n$ and $\forall i<m\, (q_i \preceq q_m \imply q_i \preceq
q_n)$. If $m>n$ then $q_m \preceq q_n$ does not hold.

It is clear that $\preceq$ is a partial ordering and its definition is
continuous (indeed uniformly primitive recursive) in $Q$ and $R$ and that
$\preceq \; \subseteq R$. If $\preceq$ is bqo then (since $\preceq \;
\subseteq R$) $R$ is bbr.

Now suppose that $\preceq$ is not bqo: since $\preceq$ is well founded it is
compatible with itself and by theorem~\ref{mba} there exists a minimal bad
$Q$-array $f$. We may assume that $\dom(f) = \Sq \N$ (this is not restrictive
since $\Sq A$ is homeomorphic to $\Sq \N$ for every $A \in \Sq \N$).

Define $h: \Sq \N \to \sq Q$ by $h(X) = \set{q \in Q}{q \prec f(X)}$: we have
$h(X) \in \sq Q$ because it has at most $\#(f(X))$ elements. We define a
partial ordering on $\sq Q$ by setting $a \preceq^* b$ if and only if
$\forall q \in a \; \exists q' \in b \; q \preceq q'$.

\begin{sublemma}
There exists $A \in \Sq \N$ such that $h$ restricted to $A$ is perfect with
respect to $\preceq^*$.
\end{sublemma}
\begin{pf}
By corollary~\ref{bad/perfect} there exists $A$ such that $h$ restricted to
$A$ is either perfect or bad. Suppose the latter holds and notice that
$\subseteq$ is well founded on $\sq Q$ and compatible with $\preceq^*$. By
theorem~\ref{mba} there exists $h' \subseteq h$ which is minimal (with
respect to $\subseteq$) bad (with respect to $\preceq^*$). Let $\dom(h') =
\Sq{A'}$: since $h'$ is bad $h'(X) \neq \emptyset$ for every $X \in \Sq{A'}$
and we can define $g(X)$ to be an element of $h'(X)$ (e.g.\ the first one to
appear in $\{q_n\}$) and $h''(X) = h'(X) \setminus \{g(X)\}$. The minimality
of $h'$ and corollary~\ref{bad/perfect} imply that $h''$ restricted to some
$A'' \in \Sq{A'}$ is perfect. Therefore $g$ restricted to $A''$ must be bad.
But for all $X \in \Sq{A''}$ we have $g(X) \in h'(X) \subseteq h(X)$ and
hence $g(X) \prec f(X)$, that is $g \prec f$ violating the minimality of $f$.
\end{pf}

Let $A$ be given by the sublemma and consider $f$ restricted to $A$:
supposing towards a contradiction that $R$ is bbr we have by
lemma~\ref{intersection} that also the intersection of $R$ with
$\set{(q,q')}{\#(q) \leq \#(q')}$ is bbr, and therefore there exists $X \in
\Sq A$ such that $f(X) \R f(\m X)$ and $\#(f(X)) \leq \#(f(\m X))$. Since
$h(X) \preceq^* h(\m X)$ we have that for every $q \prec f(X)$ there exists
$q' \prec f(\m X)$ such that $q \preceq q'$ and hence
$$\forall i < \#(f(X)) \left( q_i \preceq f(X) \imply q_i \preceq f(\m X)
\right)$$
Therefore $f(X) \preceq f(\m X)$, against the badness of $f$.
\end{pf}

\section{Smooth subsets of $\sq \N$}\label{sec:smooth}

In this section we study the notion of smooth subset of $\sq \N$, which was
introduced in \cite{thesis} and \cite{fbqo} to solve some problems in the
fine theory of bqos. Every $C \subseteq \sq \N$ can be viewed as a code for
a collection of basic open (indeed clopen) subsets of $\Sq \N$ (namely
$\{N_s\}_{s \in C}$).

\begin{definition}
If $C \subseteq \sq \N$ let $\base(C) = \set{n}{\exists s \in C \; n \in s}$.

$C$ is a {\em block\/} if $\base(C)$ is infinite, the elements of $C$ are
mutually incomparable under $\init$ and $\forall X \in \Sq{\base(C)} \;
\exists s \in C \; s \init X$.
\end{definition}

It is immediate to check that $C$ is a block if and only if $\{N_s\}_{s \in
C}$ is a partition of $\Sq{\base(C)}$. The following lemma shows the
connection between blocks and $Q$-arrays.

\begin{lemma}\label{constant}
If $f: \Sq A \to Q$ is a continuous $Q$-array there exists a block $B$ with
$\base(B) = A$ such that $f$ is constant on $N_s$ for each $s \in B$.
\end{lemma}
\begin{pf}
If $X \in \Sq A$, since $\{f(X)\}$ is open in $Q$, there exists $n > 0$ such
that $f(Y) = f(X)$ for every $Y \in N_{X[n]}$: denote by $n(X)$ the least
such $n$ and let $B = \set{X[n(X)]}{X \in \Sq A}$.
\end{pf}

\begin{definition}
If $s,t \in \sq\N$ let $t \db s$ mean that $\lh(t) = \lh(s)$ and for all $i
< \lh(s)$ $t(i) \leq s(i)$.

A set $C \subseteq \sq \N$ is {\em smooth\/} if for all $s,t \in C$ such that
$\lh(s) < \lh(t)$ there exists $i < \lh(s)$ such that $s(i) < t(i)$, i.e.\
$t[\lh(s)] \db s$ does not hold.
\end{definition}

We will use the following terminology for trees.

\begin{definition}
A set $T \subseteq \sq \N$ is a {\em tree\/} if it is closed under initial
segments, i.e.\ $s \in T$ and $t \init s$ imply $t \in T$. $T$ is {\em well
founded\/} if for every $X \in \Sq{\base(T)}$ there exists $n$ such that
$X[n] \notin T$.
\end{definition}

If we have a subset of $\sq \N$ we will turn it into a smooth subset by the
following procedure.

\begin{definition}\label{*}
If $C \subseteq \sq \N$ let $A = \base(C)$ and define
\begin{align*}
T(C) &= \set{s \in \sq A}{\forall t \initeq s \; t \notin C}\\
\T &= \set{s \in \sq A}{\exists t \in T(C)\; t \db s}\\
C^* &= \set{s \in \sq A}{s \notin \T \land \forall t \init s \; t \in \T}
\end{align*}
\end{definition}

Notice that both $T(C)$ and $\T$ are trees and $T(C) \subseteq \T$.

\begin{lemma}\label{*smooth}
$C^*$ is smooth for every $C \subseteq \sq \N$.
\end{lemma}
\begin{pf}
Assume, towards a contradiction, that there exist $s,s' \in C^*$ such that
$\lh(s) < \lh(s')$ and $s'[\lh(s)] \db s$. $s'[\lh(s') -1] \in \T$ implies
that for some $t' \in T(C)$ we have $t' \db s'[\lh(s') -1]$. Let $t =
t'[\lh(s)]$: then $t \in T(C)$ and $t \db s$ which entail $s \in \T$. This
contradicts $s \in C^*$.
\end{pf}

\begin{lemma}\label{partition}
If $C$ is a block then:
\begin{enumerate}
\item $T(C)$ and $\T$ are well founded trees.
\item For every $t \in C$ there exists $s \in C^*$ such that $t \initeq s$
and hence $N_s \subseteq N_t$.
\item $C^*$ is a block.
\end{enumerate}\end{lemma}
\begin{pf}
(1) Let $A = \base(C)$. Since $\forall X \in \Sq A \; \exists s \in C \; s
\init X$ it is obvious that $T(C)$ is well founded. Suppose $\T$ is not well
founded and let $X \in \Sq A$ be such that $\forall i\; X[i] \in \T$. Then
for every $i$ there exists $t \in T(C)$ such that $t \db X[i]$. Hence the
tree $\set{t \in T(C)}{t \db X[\lh(t)]}$ is an infinite finitely branching
tree (a node $t$ has at most $X(\lh(t)) + 1$ immediate successors) and by
K\"onig's lemma is not well founded, which is impossible because $T(C)$ is
well founded.

(2) If $t \in C$ then, since the elements of $C$ are incomparable under
$\init$, $t[\lh(t) -1] \in T(C) \subseteq \T$: by (1) there exists $s \in
C^*$ such that $t \initeq s$.

(3) By (1) $\base(C^*) = \base (C)$ and hence $\base(C^*)$ is infinite, It
is clear that two elements of $C^*$ are incomparable under $\init$. Let $X
\in \Sq A$: by (1) there exists $k$ such that $X[k] \notin \T$. If $k$ is
minimal then $X[k] \in C^*$.
\end{pf}

The binary relation we will now introduce is basic to {\sc Nash-Williams}'
combinatorial definition of bqo, but here will be employed only as a
technical tool in the proof of theorem~\ref{main lemma}.

\begin{definition}
Let $s,t \in \sq \N$: we write $s \tri t$ if there exists $u \in \sq \N$ such
that $s \initeq u$ and $t \initeq \m u$ or, equivalently, if there exists $X
\in \Sq \N$ such that $X \in N_s$ and $\m X \in N_t$.
\end{definition}

The following lemmas follow immediately from the definitions.

\begin{lemma}\label{tri}
$s \tri t$ implies $s(i) < t(i)$ and $s(i+1) = t(i)$ whenever these
expressions make sense.
\end{lemma}

\begin{lemma}\label{smooth-tri}
If $C$ is smooth and $s,t \in C$ are such that $s \tri t$ then $\lh(s) \leq
\lh(t)$.
\end{lemma}

In the proof of theorem~\ref{main lemma} we will need the following result
about $C$-arrays for smooth sets $C$ with the binary relation $\tri$.

\begin{lemma}\label{tri-pres}
Suppose $C \subseteq \sq \N$ is smooth and $f: \Sq A \to C$ is a continuous
$C$-array which is perfect with respect to $\tri$. By lemma~\ref{constant}
let $B$ be a block with $\base(B) = A$ such that $f$ has constant value
$\f(s)$ on $N_s$ for each $s \in B$. Then for every $s \in B$ we have:
\begin{enumerate}
\item $\lh(\f(s)) \leq \lh(s)$;
\item $\forall i < \lh(\f(s)) \; s(i) \leq \f(s)(i)$, i.e.\ $s[\lh(\f(s))]
\db \f(s)$.
\end{enumerate}\end{lemma}
\begin{pf}
Without loss of generality we may assume that $A = \N$. Recalling the
definition of $\tri$ notice that $f$ perfect implies $\forall s,t \in B \;
(s \tri t  \imply \f(s) \tri \f(t))$.

(1) Let $\lh(s) = k$: for every $i \leq k$ there exist unique $s_i, s'_i \in
B$ such that
\begin{align*}
s_i &\init \lt s(i), \ldots, s(k-1), s(k-1) +1, s(k-1) +2, \ldots \rt \\
s'_i &\init \lt s(i), \ldots, s(k-1), s(k-1) +2, s(k-1) +3, \ldots \rt
\end{align*}
Therefore $s = s_0 = s'_0$, $s_k \tri s'_k$ and for every $i<k$ $s_i \tri
s_{i+1}$ and $s'_i \tri s'_{i+1}$. Suppose that $\lh(\f(s)) >k$ so that
$\f(s)(k)$ exists: $C$ is smooth and by lemma~\ref{smooth-tri} (insuring that
all the sequences involved are long enough) and lemma~\ref{tri} we have
\begin{align*}
\f(s)(k) &= \f(s_1)(k-1) = \cdots = \f(s_k)(0)\\
\f(s)(k) &= \f(s'_1)(k-1) = \cdots = \f(s'_k)(0)
\end{align*}
Therefore $\f(s_k)(0) = \f(s'_k)(0)$, contradicting $\f(s_k) \tri \f(s'_k)$.

(2) We prove this simultaneously for all $s \in B$ by induction on $i$. If
$i=0$, denoting again $\lh(s)$ by $k$, for every $j \leq s(0)$ let $s_j \in
B$ be such that
$$s_j \init \lt s(0)-j, s(0)-j+1, \ldots, s(0)-1) \rt \conc s \sconc \lt
s(k-1) + 1, s(k-1) + 2, \ldots \rt$$
Therefore $s_0 = s$ and $s_{j+1} \tri s_j$ for every $j < s(0)$. By
lemma~\ref{tri} we have
$$\f(s)(0) \geq \f(s_1)(0)+1 \geq \cdots \geq \f(s_{s(0)})(0) + s(0) \geq
s(0)$$

Now suppose $i+1 < \lh(\f(s))$ and pick $t \in B$ such that $s \tri t$ (such
a $t$ exists because $B$ is a block): we have $\f(s) \tri \f(t)$ and, by
lemma~\ref{smooth-tri}, $i < \lh(\f(t))$. By lemma~\ref{tri} and the
induction hypothesis $\f(s)(i+1) = \f(t)(i) \geq t(i) = s(i+1)$.
\end{pf}

\section{Proof of the \P12-completeness of bqo}\label{sec:main}

\begin{theorem}\label{main lemma}
Let $P \subset \Bai$ be a \P12 set. For each $x \in \Bai$ we can define a
countable set $Q_x$ and a reflexive binary relation $R_x$ on $Q_x$ such that
$x \in P$ if and only if $R_x$ is a bbr. Moreover, if we view $(Q_x,R_x)$ as
an element of $\Can \times \Can \cong \Can$, the map $x \mapsto (Q_x,R_x)$
is continuous.
\end{theorem}
\begin{pf}
By the representation theorem for \S11 sets (\cite{Mosch}) and using the fact
that $\Sq \N$ is isomorphic to $\Bai$ there exists $C \subset \bigcup_n
\omega^n \times 2^n \times \left[ \N \right]^n$ such that
$$x \in P \sse \forall y \in \Can \; \exists X \in \Sq \N \; \forall n \;
(x[n], y[n], X[n])\; \notin C$$
We suppose that $(\lt \rt, \lt \rt, \lt \rt) \;\notin C$ (otherwise $P =
\emptyset$ and for every $x$ we can take $R_x$ to be a fixed non-bbr
relation). Moreover we can suppose that $C$ consists of sequences
incomparable under $\init \times \init \times \init$.

Throughout this proof $\s$ and $\tau$ will denote sequences in $2^{<\omega}$,
while $s$ and $t$ will denote elements of $\sq \N$.

Let $x \in \Bai$ be fixed and define $C_x = \set{(\s,s)}{(x[\lh(\s)], \s, s)
\in C}$ and
$$\C = \set{(\s,s)}{\lh(\s) = \lh(s) \land \forall t \db s\; \exists i \leq
\lh(s) \; (\s[i], t[i]) \in C_x)}$$
$\C$ is countable and can be viewed as an element of $\Can$: then the map $x
\mapsto \C$ is continuous because for every $s$ the set $\set{t}{t \db s}$
is finite. Let
$$Q_x = \set{(\s,s) \in \C}{\forall i < \lh(s) \, (\s[i], s[i]) \notin \C}$$
and define a binary relation $R_x$ on $Q_x$ by:
$$(\s, s) \,R_x\, (\tau, t) \sse \s \not\initeq \tau \lor s \ntri t$$
Notice that the map $x \mapsto (Q_x,R_x)$ is continuous.

To complete the proof of the theorem by showing that $x \in P$ if and only
if $R_x$ is bbr we need a couple of definitions and a sublemma.

For any $y \in \Can$ let $\CF = \set{s}{(y[\lh(s)], s) \in C_x}$. Notice that
the elements of $\CF$ are incomparable under $\init$ and, using
definition~\ref{*}, construct $\Cf$.

\begin{sublemma}
For every $y \in \Can$ and $s \in \sq\N$ we have $(y[\lh(s)],s) \in Q_x$ if
and only if $s \in \Cf$.
\end{sublemma}
\begin{pf}
If $(y[\lh(s)],s) \in Q_x$ then $(y[\lh(s)], s) \in \C$ and $\forall i <
\lh(s) \, (y[i], s[i]) \notin \C$: therefore for every $t \db s$ there exists
$i \leq \lh(s)$ such that $(y[i], t[i]) \in C_x$ and for at least one $t$
this $i$ is $\lh(s)$. Hence $\forall t \db s\; \exists i \leq \lh(s)\; t[i]
\in \CF$ and $\exists t \db s\; t \in \CF$. Therefore $\forall t \db s \; t
\notin T(\CF)$ and $\exists t \db s \; \forall i < \lh(t) \; t[i] \in
T(\CF)$. Hence $s \notin T^*(\CF)$ and $\forall i < \lh(s) \; s[i] \in
T^*(\CF)$: thus $s \in \Cf$.

If $s \in \Cf$ we have $s \notin T^*(\CF)$. Hence if $t \db s$ we have $t
\notin T(\CF)$ and there exists $i \leq \lh(t)$ such that $t[i] \in \CF$,
i.e.\ $(y[i], t[i]) \in C_x$: this means $(y[\lh(s)],s) \in \C$. On the other
hand for all $i < \lh(s)$ we have $s[i] \in T^*(\CF)$ and hence there exists
$t \db s[i]$ such that $t \in T(\CF)$: for every $j \leq i$ we have $t[j]
\notin \CF$ and hence $(y[j], t[j]) \notin C_x$. This shows that $(y[i],
s[i]) \notin \C$ and completes the proof that $(y[\lh(s)],s) \in Q_x$.
\end{pf}

First assume that $x \notin P$: for some $y \in \Can$ we have $\forall X \in
\Sq \N\; \exists n\; X[n] \in \CF$. Hence $\CF$ is a block and, by
lemma~\ref{partition}.3, $\Cf$ is also a block. If $X \in \Sq \N$ denote by
$s_X$ the unique element of $\Cf$ such that $s_X \init X$, so that $s_X \tri
s_{\m X}$. Define $f: \Sq \N \to Q_x$ by $f(X) = \left( y[\lh(s_X)], s_X
\right)$. By the sublemma $f(X) \in Q_x$ and it is immediate, using
lemma~\ref{smooth-tri}, to check that $f$ is a bad continuous $Q_x$-array
with respect to $R_x$ and hence $R_x$ is not bbr.

Suppose $x \in P$ but $R_x$ is not bbr: by corollary~\ref{cont=borel} there
exists a continuous bad $Q_x$-array $f: \Sq \N \to Q_x$. Let $B$ and $\f$ be
as in the statement of lemma~\ref{tri-pres} so that $\f: B \to Q_x$. Let us
write $\f(s) = (g(s), h(s))$. The badness of $f$ implies that whenever $s,
t \in B$ are such that $s \tri t$ we have (a)~$g(s) \initeq g(t)$ and
(b)~$h(s) \tri h(t)$.

If $\tri^*$ is the transitive closure of $\tri$ and $s,t \in B$ there exists
$u \in B$ such that $s \tri^* u$ and $t \tri^* u$. By~(a) we have $g(s)
\initeq g(u)$ and $g(t) \initeq g(u)$ so that either $g(s) \initeq g(t)$ or
$g(t) \initeq g(s)$. It follows that there exists $y \in \Can$ such that
$\forall s \in B \; g(s) \init y$ and hence $\f(s) = (y[\lh(h(s))], h(s)) \in
Q_x$.

By the sublemma the range of $h$ is a subset of $\Cf$. Since $x \in P$ there
exists $X \in \Sq \N$ such that $\forall n \; (y[n], X[n]) \notin C_x$. If
$s \in B$ is such that $s \init X$, by~(b) and since $\Cf$ is smooth by
lemma~\ref{*smooth}, lemma~\ref{tri-pres}.2 applies and $X[\lh(h(s))] \db
h(s)$: it is then obvious that $f(X) = \left( y[\lh(h(s))], h(s) \right)
\notin \C$, a contradiction.
\end{pf}

The following theorem is the main result of the paper.

\begin{theorem}\label{main theorem}
Let $BQO$ be the set of all (codes for) bqos on $\N$. $BQO$ is a
\P12-complete subset of $\Can$.
\end{theorem}
\begin{pf}
By lemma~\ref{constant} a $\N$-array $f: \Sq A \to \N$ is continuous if and
only if there exists a block $B$ with $\base(B) = A$ such that $f$ is
constant on $N_s$ for every $s \in B$: the set of all blocks is a \P11 subset
of $\Can$ and using corollary~\ref{cont=borel} it is immediate to see that
$BQO$ is a \P12 subset of $\Can$.

To show that $BQO$ is \P12-complete, given a \P12 set $P \subseteq \Bai$ it
suffices to apply theorem~\ref{main lemma} followed by
theorem~\ref{th:pouzet} to obtain a reduction of $P$ to $BQO$.
\end{pf}

\end{document}